%% file: rank.tex
\documentclass[12pt,a4paper]{article}

\usepackage{amssymb,amsfonts,amsmath,amsthm,txfonts}

\usepackage{hyperref}

\newtheorem{theorem}{Theorem}
\newtheorem{lemma}[theorem]{Lemma}
\newtheorem{proposition}[theorem]{Proposition}
\newtheorem{corollary}[theorem]{Corollary}

\theoremstyle{definition}
\newtheorem{definition}{Definition}
\newtheorem{example}{Example}

\newcommand{\nrank}{\ensuremath{N\text{-rank}}}
\newcommand{\ord}{\ensuremath{\text{ord}}}

\title{A lower bound for the $r$-order of a matrix modulo N\footnotetext{2000 Mathematics Subject Classification: 11C20, 11G20, 14G05.}\footnotetext{Key words: matrices, rank, order, elliptic curves, finite fields, rational points.}}
\author{Carlo Magagna}
\date{}

\begin{document}
\maketitle

\input{content}

\bibliographystyle{abbrv}
\bibliography{rank}

\noindent Carlo Magagna\\
Dipartimento di Matematica e Informatica\\
Via Delle Scienze 206\\
33100, Udine, Italy\\
email: magagna@dimi.uniud.it

\end{document}

%% file: content.tex
\begin{abstract}
For a positive integer $N$, we define the $\nrank$ of a non singular integer $d\times d$ matrix $A$ to be the maximum integer $r$ such that there exists a minor of order $r$ whose determinant is not divisible by $N$. Given a positive integer $r$, we study the growth of the minimum integer $k$, such that $A^k-I$ has $\nrank$ at most $r$, as a function of $N$. We show that this integer $k$ goes to infinity faster than $\log N$ if and only if for every eigenvalue $\lambda$ which is not a root of unity, the sum of the dimensions of the eigenspaces relative to eigenvalues which are multiplicatively dependent with $\lambda$ and are not roots of unity, plus the dimensions of the eigenspaces relative to eigenvalues which are roots of unity, does not exceed $d-r-1$. This result will be applied to recover a recent theorem of Luca and Shparlinski \cite{ls:ell} which states that the group of rational points of an ordinary elliptic curve $E$ over a finite field with $q^n$ elements is almost cyclic, in a sense to be defined, when $n$ goes to infinity. We will also extend this result to the product of two elliptic curves over a finite field and show that the orders of the groups of $\mathbb{F}_{q^n}-$rational points of two non isogenous elliptic curves are almost coprime when $n$ approaches infinity.
\end{abstract}
\section{Introduction}
In \cite{cz:primefactor} it is shown that if $S$ is a finite set of absolute values of $\mathbb{Q}$, including $\infty$, and $u,v\in\mathbb{Z}$ are multiplicatively independent $S$-units, then for every $\epsilon >0$, $\gcd(u-1,v-1)<\max(|u|,|v|)^\epsilon$ holds with finitely many exceptions (see also \cite{b:heights}, theorem 7.4.10).
A more general result is presented in \cite{cz:lowheight} where an upper bound
for the greatest common divisor, to be defined later, between $u-1$ and $v-1$ where $u,v$ are now $S$-units in a number field, is obtained. This result is used in \cite{cz:matperiod} to show that for a non singular integer matrix $A$, the growth of the order of $A$ modulo an integer $N$ goes to infinity faster then $\log N$ if and only if none of the following cases holds:
\begin{itemize}
\item[(i)] $A$ is diagonalizable and a power of $A$ has all the eigenvalues equal to powers of a single rational integer
\item[(ii)] $A$ is diagonalizable and a power of $A$ has all the eigenvalues equal to powers of a single unit in a real quadratic field
\end{itemize}

Let now $d$ be a positive integer and $A$ a non singular $d \times d$ integer
matrix. Given an integer $N\ge 1$, we define the $\nrank$ of $A$ as
follows.
\begin{definition}
\label{def:nrank}
The $\nrank$ of $A$ is the greatest integer $r\ge 0$ such that there exists an $r\times r $ minor of
$A$ whose determinant is not divisible by $N$. We will write
$r=\nrank(A)$.
\end{definition}
Given $r$ as in definition \ref{def:nrank} we can define the $r$-order of
the matrix $A$ as follows.
\begin{definition}
A positive integer $k$ is called the $r$-order of $A$ modulo $N$, if
it is the smallest integer such that $\nrank(A^k-I)\le r$, where $I$
denotes the identity matrix. We will write $k=\ord(A,N,r)$. If such an
integer does not exist, we will set $\ord(A,N,r)=\infty$.
\end{definition}
Before stating the results of this paper, we briefly analyze the main properties of $\nrank$ and $r$-order. An integer matrix $A$ has $\nrank$ zero if and only if all the determinants of order $1$ are divisible by $N$, i.e. if and only if $A\equiv 0 \pmod{N}$. Hence $\ord(A,N,0)$ is just the usual order $\ord(A,N)$ of a matrix and have been studied in \cite{cz:matperiod} as recalled above. Recall now that from an integer $d\times d$ matrix $A=(a_{ij})$ and a given positive integer $r\le d$, one can construct a new matrix, the so called $r$-th exterior power of $A$ as follows. Let $S_r^d$ be the set of sequences
\begin{equation*}
 J=(j_1,j_2,\ldots ,j_r),\quad \text{where }1\le j_1<j_2<\ldots <j_r\le d
\end{equation*}
and let $J,K\in S_r^d$. Then we define
\begin{equation}
\label{eq:extpow}
 A_{J,K}^{(r)}\coloneqq
\det\left( \begin{array}{cccc}
a_{j_1 k_1} & a_{j_1 k_2} & \cdots & a_{j_1 k_r}\\
a_{j_2 k_1} & a_{j_2 k_2} & \cdots & a_{j_2 k_r}\\
\vdots & \vdots &  & \vdots \\
a_{j_r k_1} & a_{j_r k_2} & \cdots & a_{j_r k_r}\\
\end{array} \right).
\end{equation}
Hence, by varying $J,K$ in $S_r^d$ we obtain a new matrix $A^{(r)}$, the $r$-th exterior power of $A$, whose $JK$-component is defined by (\ref{eq:extpow}).  Choose now a $\mathbb{Z}$-basis $\{e_1,e_2,\ldots,e_d\}$ for $\mathbb{Z}^d$; then the elements
\begin{equation*}
 e_{j_1}\wedge e_{j_2}\wedge \ldots \wedge e_{j_r},\quad \text{where }1\le j_1<j_2<\ldots <j_r\le d
\end{equation*}
form a basis for the Grassmann algebra $\bigwedge^r\mathbb{Z}^d$ and the matrix $A^{(r)}$ represents an endomorphism \begin{equation*}
{\textstyle \bigwedge^r A:\bigwedge^r\mathbb{Z}^d \longrightarrow \bigwedge^r\mathbb{Z}^d.}
\end{equation*}
Now $\nrank(A)\le r$ if and only if every determinant of a minor of order $r+1$ is divisible by $N$ and this in turn is equivalent to the fact that the matrix $A^{(r+1)}$, representing $\bigwedge^{r+1}A$, has all the entries divisible by $N$, i.e. $N|A^{(r+1)}$.

Moreover, if one defines, for $0\le r\le d$, the determinant ideal $I_r(A)$ to be the ideal in $\mathbb{Z}$ generated by the entries of $A^{(r)}$, i.e. by the determinants of the minors of $A$ of order $r$, then, after putting $I_0(A)\coloneqq\mathbb{Z}$,
\begin{equation*}
 I_0(A)\supset I_1(A)\supset \cdots\supset I_d(A).
\end{equation*}
It follows that, for every $r$ between $0$ and $d-1$,
\begin{equation}
\label{eq:monotone}
 \ord(A,N,r)\ge \ord(A,N,r+1).
\end{equation}
Let now $P,Q\in\text{GL}_d(\mathbb{Z})$. In this case $I_r(P)=I_r(Q)=(1)$ for every $r$, hence, since for two $d\times d$ integer matrices $A,B$,
\begin{equation*}
 (AB)^{(r)}=A^{(r)}B^{(r)},
\end{equation*}
then
\begin{equation*}
 I_r(A)=I_r(PAQ)
\end{equation*}
and we deduce that the $\nrank$ is invariant under conjugation in $\text{GL}_d(\mathbb{Z})$.

\section{Statements of the results}
In this paper we study the minimal growth of $\ord(A,N,r)$, for fixed
values of $r$ and given $A$, as $N\rightarrow \infty$.
If $A$ has finite $r$-order (globally), i.e. $A^k-I$ has rank at most
$r$ for a certain $k\ge 1$, then clearly $\ord(A,N,r)\le k$ is
bounded. If this is not the case, then $\ord(A,N,r)\rightarrow \infty$
as $N\rightarrow \infty$.\\
The case $r=d$ is trivial, being $\ord(A,N,d)=1$ for each integer
$N\ge 1$. When $r=d-1$ the growth is not faster then logarithmic. Let us first consider the case where no eigenvalue is a root of unity. Let $\lambda_1,\ldots,\lambda_t,\lambda_{t+1},\ldots\lambda_d$ be the complex
eigenvalues of $A$, taken with multiplicity 1 and ordered in a way
that $|\lambda_i|>1$ if and only if $i\le t$. Let
$N_n=\left|\det(A^n-I)\right|$ and
$\eta=\sum_{i=1}^{t}\log|\lambda_i|$. Observing that
$$\log N_n=\sum_{i=1}^{t}\log\left|\lambda_i^n-1\right|+O(1)=n\eta + O(1)$$
we obtain
$$\ord(A,N_n,d-1)\le n =\eta^{-1}\log N_n + O(1)$$
and so
$$\liminf_{N\rightarrow\infty} \frac{\ord(A,N,d-1)}{\log N}\le\eta^{-1}<\infty$$
as wanted.\\
If an eigenvalue, say $\lambda_1$, is a root of unity: if $\lambda_1^m=1$ then $\det(A^m-I)=0$, therefore $\ord(A,N,d-1)\le m$ for every positive integer $N$.

\vskip\baselineskip
From now on we will then consider $0\le r \le d-2$. Let $\mathbb{K}\subset\overline{\mathbb{Q}}$ be the splitting field of the characteristic polynomial of
$A$. Then there exists an invertible matrix $P$ over $\mathbb{K}$ such that
\begin{equation}
\label{eq:Jordan}
B=P^{-1} A P
\end{equation}
is the Jordan canonical form of $A$.
Let now $\Lambda$ be the set of
eigenvalues of $A$, let $\Lambda^\ast\subset\Lambda$ be the set of
eigenvalues that are roots of unity and let $\Lambda'\coloneqq\Lambda\setminus\Lambda^\ast$.
\begin{definition}
Two eigenvalues $\lambda_1,\lambda_2\in \Lambda$ are multiplicatively dependent if and only if there exists $(a_1,a_2)\in\mathbb{Z}^2\setminus\{(0,0)\}$ such that $\lambda_1^{a_1}\lambda_2^{a_2}=1$.
\end{definition}
Let $\sim$ be the equivalence
relation of being pairwise multiplicatively dependent, defined on the
set $\Lambda'$ of eigenvalues of $A$ which are not roots of unity and
let $\Gamma\coloneqq\Lambda'/\sim$. Note that $\sim$ would not be an equivalence relation if defined on the whole $\Lambda$, since every eigenvalue is multiplicatively dependent with an eigenvalue in $\Lambda^\ast$ and transitivity would fail. For each equivalence class
$\gamma\in \Gamma$ we set $h_\gamma$ to be the sum of the
algebraic multiplicities of the eigenvalues in $\gamma$ and
$\overline{h}_\gamma$ to be the number of $1$ appearing in the
Jordan blocks of $B$ relative to the eigenvalues in $\gamma$.
Finally let $l$ be the sum of the algebraic multiplicities of the
eigenvalues in $\Lambda^\ast$ and $\overline{l}$ be the number of $1$ appearing in the
Jordan blocks of $B$ relative to the eigenvalues in $\Lambda^\ast$.
\begin{definition}
Given an integer $r$, with $0\le r \le d-2$, a $d\times d$ integer
matrix $A$ will be called $r$-regular if 
\begin{equation}
\label{eq:lim}
\lim_{N\rightarrow\infty}\frac{\ord(A,N,r)}{\log N}=+\infty
\end{equation}
and $r$-exceptional otherwise.
\end{definition}
The main result of this note is the following theorem.
\begin{theorem}
\label{mainprop}
Let $A$ be a non singular integer $d\times d$ matrix and $r\le d-2$ a non negative
integer. Then $A$ is $r$-exceptional if and only if there exists
$\gamma\in \Gamma$ such that
\begin{equation}
\label{eq:mainprop}
l-\overline{l}+h_\gamma-\overline{h}_\gamma\ge d-r
\end{equation}
\end{theorem}
\begin{example}\label{ex:1}
Consider for example the matrix
\begin{equation*}
 A=\left(
\begin{array}{lllll}
 21 & -10 & 2 & -12 & 1 \\
 15 & -7 & 5 & -15 & 3 \\
 3 & -2 & 4 & -3 & 1 \\
 9 & -4 & -1 & 0 & -1 \\
 -2 & 1 & 1 & 2 & 2
\end{array}
\right).
\end{equation*}
This matrix has Jordan canonical form
\begin{equation*}
  B=\left(
\begin{array}{lllll}
 2 & 0 & 0 & 0 & 0 \\
 0 & 3 & 1 & 0 & 0 \\
 0 & 0 & 3 & 1 & 0 \\
 0 & 0 & 0 & 3 & 0 \\
 0 & 0 & 0 & 0 & 9
\end{array}
\right).
\end{equation*}
which has $2, 3$ and $9$ as eigenvalues and none of them is a root of unity. Two of them, $3$ and $9$ are multiplicatively dependent, while $2$ is multiplicatively independent with $3$ and hence with $9$. Therefore $l=\overline{l}=0$ and, if we denote with $\overline{2}$, $\overline{3}$ and $\overline{9}$ the classes in $\Gamma$ that contain respectively $2, 3$ and $9$, then $\overline{3}=\overline{9}$ and $h_{\overline{2}}=1$,  $\overline{h}_{\overline{2}}=0$, $h_{\overline{3}}=4$,
$\overline{h}_{\overline{3}}=2$. Then, by applying Theorem \ref{mainprop}, the matrix $A$ is $3$-exceptional, $2$-regular and then $1$- and $0$-regular, by equation (\ref{eq:monotone}).
\end{example}
\begin{example}Consider now an invertible integer matrix whose Jordan form is
\begin{equation*}
 A=\left( \begin{array}{ccccc}
\zeta_1 & 0 & 0 & 0 & 0\\
0 & \zeta_2 & 0 & 0 & 0\\
0 & 0 & a & 1 & 0\\
0 & 0 & 0 & a & 0 \\
0 & 0 & 0 & 0 & b \\
\end{array} \right),
\end{equation*}
where $\zeta_1,\zeta_2$ are roots of unity and $a,b$ are two multiplicatively dependent non roots of unity. In this case $l=2$, $\overline{l}=0$ and, in the notation of example \ref{ex:1}, $h_{\overline{a}}=3$,  $\overline{h}_{\overline{a}}=1$. Then, $l-\overline{l}+h_{\overline{a}}-\overline{h}_{\overline{a}}=4$, hence by applying Theorem \ref{mainprop}, the matrix $A$ is $0$-regular and $1$-exceptional (and then $2$- and $3$-exceptional, since (\ref{eq:monotone}) holds).
\end{example}

The main tool to prove the necessity of condition (\ref{eq:mainprop}) for $A$ being $r-$exceptional will be a result of diophantine approximation by Corvaja and Zannier \cite{cz:lowheight} which is an application of Schmidt's subspace theorem. On the other hand, to prove the sufficiency of (\ref{eq:mainprop}), a generalized version of Roth's theorem will suffice.\\
As an immediate corollary of Theorem \ref{mainprop} we can deduce a sufficient condition on the structure of the Zariski closure $G_A\coloneqq \overline{<A>}$ in $GL_d$ of the cyclic group generated by a single invertible integer matrix $A$, for $A$ being $r$-regular. Let $G_A^0$ be the connected component of $G_A$ containing the identity; then, by the general theory of commutative algebraic groups, $G_A^0\cong \mathbb{G}_m^e\times \mathbb{G}_a^f$, where $\mathbb{G}_m$ and $\mathbb{G}_a$ denote respectively the multiplicative and the additive groups and $f=0$ or $1$ depending on $A$ being diagonalizable or not.
\begin{corollary}
\label{cor:alggrp}
Let $A$ be an invertible integer matrix, $r\le d-2$ a non negative
integer and $G_A^0\cong\mathbb{G}_m^e\times \mathbb{G}_a^f$ the connected component containing the identity of the Zariski closure of the group generated by $A$. If $e+f>r+1$, then $A$ is $r$-regular.
\end{corollary}
The converse of the corollary is not true. Consider for example a $3\times 3$ diagonalizable matrix $A$ with three distinct eigenvalues $\lambda,\mu,\nu$ non multiplicative dependent in pairs, but such that there exist three integers $a,b,c$ such that $\lambda^a\mu^b\nu^c=1$, for instance $\lambda=3,\mu=5,\nu=15$. Then $A$ is $1$-regular, but $e+f=2$.

\vskip\baselineskip
For certain applications it is more convenient to consider, more generally then an integer matrix, an endomorphism $\phi$ of a finitely generated free module over a ring of characteristic zero, without choosing a base. If all the coefficients of the characteristic polynomial of $\phi$ are rational integers, then such are the coefficients of the characteristic polynomial of $\phi^n-I$, for every positive integer $n$, where $I$ is the identity endomorphism. In the following the coefficients of the characteristic polynomial of an endomorphism $\phi$, will be called the {\em invariants} of $\phi$. Let us, for every $k=1,\ldots,d$, denote with $\alpha_{n,k}$ the invariant of $\phi^n-I$ which is homogeneous of degree $k$ in the eigenvalues of $\phi^n-I$. We can then consider, for fixed $N\in\mathbb{N}$, the smallest positive integer $k(\phi,N)$ such that $N$ divides $\alpha_{n,k}^{d!k^{-1}}$ for all $k=1,\ldots,d$. A slight modification of the arguments used in proving Theorem \ref{mainprop} leads to the following result.
\begin{theorem}
\label{invarprop}
Let $\phi$ be an endomorphism of a finitely generated free module over a ring of characteristic zero, such that the invariants of $\phi$ are rational integers. Then $k(\phi,N)$, defined as above, satisfies
\begin{equation}
\label{eq:liminv}
\lim_{N\rightarrow\infty}\frac{k(\phi,N)}{\log N}=+\infty
\end{equation}
if and only if $\phi$ has at least two multiplicatively independent eigenvalues. 
\end{theorem}
As an application of this theorem we can recover a result of Luca and Shparlinski, presented in \cite{ls:ell}, on the exponent of the group of rational points on an elliptic curve defined over a finite field. Let $E$ be an elliptic curve defined over a finite field $\mathbb{F}_q$, with $q$ elements, and let $E(\mathbb{F}_{q^n})$ be the group of $\mathbb{F}_{q^n}-$rational points. It is known that $E(\mathbb{F}_{q^n})$ has the following structure \cite[chapter 5]{silver:arit_ell}:
\begin{equation}
\label{eq:grpstruct}
E(\mathbb{F}_{q^n})\cong \left(\mathbb{Z}/m(q^n)\mathbb{Z}\right) \times \left(\mathbb{Z}/l(q^n)\mathbb{Z}\right)
\end{equation}
where $m(q^n), l(q^n)$ are uniquely determined integers such that $m(q^n)|l(q^n)$. The integer $l(q^n)$ is the largest possible order of torsion of an $\mathbb{F}_{q^n}-$rational point and it is called the {\em exponent} of $E(\mathbb{F}_{q^n})$. Moreover the Hasse-Weil relation for the cardinality $\sharp E(\mathbb{F}_{q^n})$ of the set of $\mathbb{F}_{q^n}-$rational points is
\begin{equation}
\label{eq:hasseweil}
\sharp E(\mathbb{F}_{q^n})=q^n+1-\text{Tr}(\phi^n)
\end{equation}
where $\phi$ is the Frobenius isogeny of $E$ and $\text{Tr}(\phi^n)$ is the trace of its $n$-th power. Using equation (\ref{eq:hasseweil}) and the fact that the eigenvalues $\alpha, \beta$ of $\phi$ are complex coniugates with $|\alpha|=|\beta|=q^{1/2}$, it is immediate to obtain the bound
\begin{equation}
\label{eq:trivialbound}
l(q^n)\ge q^{n/2}-1
\end{equation}
for every $n$. We will apply Theorem \ref{invarprop} to recover the much stronger lower bound of Luca and Shparlinski for the exponent of $E(\mathbb{F}_{q^n})$ for an ordinary elliptic curve. To state their theorem, recall that an elliptic curve defined over $\mathbb{F}_q$, with $q=p^k$, is said ordinary if the group of $p$-torsion points is isomorphic to $\mathbb{Z}/p\mathbb{Z}$ and supersingular if $0$ is the unique $p$-torsion point.
\begin{theorem}[F. Luca and E. Shparlinski]
\label{ellprop}
Let $E$ be an elliptic curve over a finite field $\mathbb{F}_q$. Then for every $\epsilon >0$,
$$l(q^n)\ge q^{n(1-\epsilon)}\;\;
\text{for every }n\text{ sufficiently large}$$
if and only if $E$ is ordinary.
\end{theorem}
To prove Theorem \ref{ellprop}, we will apply Theorem \ref{invarprop} with $\phi$ equals to the Frobenius endomorphism of the elliptic curve $E$, by showing that $m(q^n)$, in the notation of (\ref{eq:grpstruct}), divides both $\det(\phi^n-I)$ and $(\text{Tr}(\phi^n-I))^2$. This fact, together with the Hasse-Weil relation (\ref{eq:hasseweil}), will provide the desired result.

\vskip\baselineskip
Consider now two ordinary elliptic curves $E_1$ and $E_2$ defined over $\mathbb{F}_{q}$. Let $\mathcal{A}\coloneqq E_1\times E_2$ be their product and let $\mathcal{A}(\mathbb{F}_{q^n})=E_1(\mathbb{F}_{q^n})\times E_2(\mathbb{F}_{q^n})$ be the group of its $\mathbb{F}_{q^n}-$rational points. Since
\begin{equation*}
E_i(\mathbb{F}_{q^n})\cong \left(\mathbb{Z}/ m_i(q^n)\mathbb{Z}\right) \times \left(\mathbb{Z}/l_i(q^n)\mathbb{Z}\right)
\end{equation*}
for $i=1,2$, then
\begin{equation*}
\mathcal{A}(\mathbb{F}_{q^n})\cong \left(\mathbb{Z}/m_1(q^n)\mathbb{Z}\right) \times \left(\mathbb{Z}/l_1(q^n)\mathbb{Z}\right)\times
\left(\mathbb{Z}/m_2(q^n)\mathbb{Z}\right) \times \left(\mathbb{Z}/l_2(q^n)\mathbb{Z}\right)
\end{equation*}
and then
\begin{equation*}
\mathcal{A}(\mathbb{F}_{q^n})\cong \left(\mathbb{Z}/l(q^n)\mathbb{Z}\right) \times M(q^n)
\end{equation*}
where
\begin{equation*}
l(q^n)\coloneqq \text{lcm}(l_1(q^n),l_2(q^n))
\end{equation*}
is the least common multiple of the exponents of the groups $E_1(\mathbb{F}_{q^n})$ and $E_2(\mathbb{F}_{q^n})$ and $M(q^n)$ is a finite, not necessarily cyclic, group.
We will apply Theorems \ref{mainprop} and \ref{ellprop} to prove the following necessary and sufficient condition on the structure of $\mathcal{A}(\mathbb{F}_{q^n})$ for the two curves to be isogenous.
\begin{theorem}
\label{2ellprop}
Let $E_1$ and $E_2$ be two ordinary elliptic curves over a finite field $\mathbb{F}_q$. Then for every $\epsilon >0$,
\begin{equation}
\label{eq:2cycl}
l(q^n)\ge q^{2n(1-\epsilon)}\exp(-\epsilon n)\;\;
\text{for every }n\text{ sufficiently large}
\end{equation}
if and only if $E_1$ and $E_2$ are not isogenous over $\overline{\mathbb{F}}_q$. Hence if $E_1$ and $E_2$ are not isogenous then
\begin{equation}
\label{eq:gcdnoiso}
\gcd \left(\sharp E_1(\mathbb{F}_{q^n}),\sharp E_2(\mathbb{F}_{q^n})\right)<\exp(\epsilon n)\;\; \text{for every }n\text{ sufficiently large.}
\end{equation}
\end{theorem}
Equation (\ref{eq:gcdnoiso}) can be paraphrased by saying that the groups of $\mathbb{F}_{q^n}-$rational points of two ordinary non isogenous elliptic curves have orders which tend to be coprime as $n$ approaches infinity.
\subsection*{Acnowledgements}
This article will be part of my PhD thesis, supervised by prof. Umberto Zannier. I would like to thank him for his helpful advice. I would like also to thank prof. Pietro Corvaja for the many useful discussions.
\section{Proofs}
To prove Theorem \ref{mainprop}, we need a few lemmas.
\begin{lemma}
Let $A$ and $r$ be as in Theorem \ref{mainprop}, $A$ not of finite
global $r$-order, $n$ a positive
integer and let $x_{n,r,i}$, $i=1,\ldots, \binom{d}{r}^2$ be the
determinants of the minors of $A^n-I$ of order $r$. Then the
following statement is equivalent to (\ref{eq:lim})
\begin{equation}
\label{eq:gcd}
\forall \epsilon > 0,\;\; \gcd_{i}(x_{n,r+1,i})< \exp (\epsilon n)\;\;
\text{for }n\text{ sufficiently large with respect to $\epsilon$.}
\end{equation}
\end{lemma}
\begin{proof}
Let $k\coloneqq\ord(A,N,r)$. Then $N$ divides $x_{k,r+1,i}$ for every $i$. In particular
\begin{equation}
\label{eq:Ngcd}
N\le \gcd_i(x_{k,r+1,i}),\quad\forall N\in\mathbb{N}
\end{equation}
If condition (\ref{eq:gcd}) holds, then
$$\gcd_i(x_{k,r+1,i})< \exp(\epsilon k),\quad\text{for $N$ (and thus
$k$) sufficiently large}.$$
Combining this with (\ref{eq:Ngcd}) we obtain
$$\frac{\ord(A,N,r)}{\log N} > \epsilon^{-1},\quad\text{for $N$
sufficiently large}$$
and this implies condition (\ref{eq:lim}).\\
On the other hand if there exist a positive real number $\rho$ and an infinite subset
$\mathcal{N}$ of $\mathbb{N}$ such that
$$\gcd_i(x_{n,r+1,i})\ge \exp(\rho n),\quad\forall n\in\mathcal{N}$$
then, taking $N_n\coloneqq \gcd_i(x_{n,r+1,i})$, we get
$$\ord(A,N_n,r)\le n \le \frac{1}{\rho}\log
\gcd_i(x_{n,r+1,i})=\frac{1}{\rho}\log N_n$$
and so
$$\frac{\ord(A,N_n,r)}{\log N_n}\le \frac{1}{\rho},\quad\forall n\in\mathcal{N}.$$
\end{proof}
We need now to introduce some notation related with $\mathbb{K}$, the splitting field of the charachteristic polynomial of $A$. Let $M$ and $M_0$ be respectively the set of places and
finite places of the field $\mathbb{K}$ and normalize the associated absolute
values in such a way that the product formula $\prod_{\mu\in M}|x|_\mu
=1$ holds for each $x\in\mathbb{K}^*$. We will also need the absolute
logarithmic Weil height $h(x)=\log H(x)$ of a point $x\in\mathbb{K}$, where
$H(x)\coloneqq \prod_{\mu\in M}\max\{1,|x|_\mu\}$. If $\{x_1,\ldots,x_k\}\subset \mathcal{O}_\mathbb{K}$ is a finite set of algebraic integers of $\mathbb{K}$, we define
$$\log\gcd_{i}(x_i)\coloneqq
\sum_{\mu\in M_0}\log^- \max_{i}\{|x_i|_\mu\}$$
to extend the concept of $\gcd$ from the rational integers to the
ring $\mathcal{O}_\mathbb{K}$ of algebraic integers of $\mathbb{K}$ and
$\log^-(x)\coloneqq -\min\{0,\log(x)\}$ for every $x>0$. Finally let $S$
be a finite subset of $M$, including $M\setminus M_0$, and let
$$\mathcal{O}_{\mathbb{K},S}^*=\{x\in\mathbb{K} \text{ such that } |x|_\mu=1,  \forall
\mu\notin S\}$$
be the group of $S$-units of $\mathbb{K}$.

\vskip\baselineskip
Noting that (\ref{eq:Jordan}) implies $A^n-I=P(B^n-I)P^{-1}$ and letting $y_{n,r,i}$,
$i=1,\ldots, \binom{d}{r}^2$ be the determinants of the minors of
$B^n-I$ of order $r$, we observe that condition (\ref{eq:gcd}) (and
thus condition (\ref{eq:lim})) holds if and only if a similar condition holds for the matrix $B$, i.e. (\ref{eq:gcd}) is equivalent to
\begin{equation}
\label{eq:lgcd}
\forall \epsilon > 0,\;\; \log\gcd_{i}(y_{n,r+1,i})< \epsilon n\;\;
\text{for }n\text{ sufficiently large.}
\end{equation}
To prove the equivalence of (\ref{eq:gcd}) and (\ref{eq:lgcd}) observe that the
entries of $P$ are fixed, independently of the exponent $n$, and hence
have bounded denominators as $n$ varies. So for
each $i=1,\ldots, \binom{d}{r}^2$, $y_{n,r,i}$ is a linear
combination of the $x_{n,r,j}$, $j=1,\ldots, \binom{d}{r}^2$ with
coefficients having bounded denominators and so $|y_{n,r,i}|_\mu \le
c_\mu \max_j |x_{n,r,j}|_\mu$, where $c_\mu=1$ for all but finitely
many $\mu\in M$. This implies the equivalence of (\ref{eq:gcd}) and
(\ref{eq:lgcd}).\\
To prove Theorem \ref{mainprop}, we begin by considering the special case of two multiplicatively dependent eigenvalues. In this case we can prove the following lemma, whose proof is elementary, in the sense that, it does not use any tool of diophantine approximation.
\begin{lemma}
\label{lem:ckgen}
Let $\lambda,\eta\in\mathbb{K}^\times$ multiplicatively dependent algebraic integers, $\lambda$ being not a root of unity, and $B(\eta)$ be a
Jordan block of order $k+1$ with exactly $k$ ``$1$'' off-diagonal:
$$B(\eta)\coloneqq
\left( \begin{array}{ccccc}
\eta & 1 & 0 & \cdots & 0\\
0 & \eta & 1& \cdots & 0\\
\vdots & \vdots & \ddots &\ddots & \vdots \\
0& \cdots & 0 & \eta & 1 \\
0 & \cdots & \cdots & 0 & \eta
\end{array} \right)$$
Let $C_{n,k}(\eta)$ be the $k\times k$ minor of $B(\eta)^n-I$ made up with
the first $k$ rows and columns $2,3,\ldots,k$. Then
\begin{equation*}
\log\gcd(\lambda^n-1, \det C_{n,k}(\eta))=O(\log n).
\end{equation*}
\end{lemma}
\begin{proof}
{\em Case 1)}
Consider first the case where $\eta$ is not a root of unity. Let $a,b$ be non zero integers such that $\lambda^a=\eta^b$. If $ab<0$, then $\lambda$ is a unity and, since $\lambda^n-1=-\lambda^n(\lambda^{-n}-1)$, the ideals generated by $\lambda^n-1$ and $\lambda^{-n}-1$ coincide, hence
$$\log\gcd(\lambda^n-1, \det C_{n,k}(\eta))=\log\gcd(\lambda^{-n}-1, \det C_{n,k}(\eta))$$
We can therefore suppose $a$ and $b$ positive, by replacing $\lambda$ with $\lambda^{-1}$ if necessary. There exists then an algebraic integer $\xi\in\mathbb{K}\left[\sqrt[b]{\lambda}\right]$ such that $\xi^b=\lambda$ and $\xi^a=\eta$. If we now set $t=\xi^n$ we get
$$\lambda^n-1=\xi^{bn}-1=t^b-1$$
and
$$\det C_{n,k}(\eta)=\eta^{-k}\det\left( \begin{array}{lllll}
n t^a & \dbinom{n}{2} t^a & \cdots & \cdots & \dbinom{n}{k} t^a\\
t^a-1 & n t^a & \cdots& \cdots & \dbinom{n}{k-1} t^a\\
\vdots & \ddots & \ddots &\vdots & \vdots \\
0& \cdots & t^a-1 & n t^a &  \dbinom{n}{2} t^a\\
0 & \cdots & \cdots & t^a-1 & n t^a
\end{array} \right)$$
It is now convenient to define two polynomials  $f,g\in\mathbb{Q}\left[x,t\right]$, with $x$ and $t$ algebraically independent over $\mathbb{Q}$, as follows:
$$f(x,t)\coloneqq t^b-1$$
$$g(x,t)\coloneqq 
\det\left( \begin{array}{lllll}
x t^a & \dbinom{x}{2} t^a & \cdots & \cdots & \dbinom{x}{k} t^a\\
t^a-1 & x t^a & \cdots& \cdots & \dbinom{x}{k-1} t^a\\
\vdots & \ddots & \ddots &\vdots & \vdots \\
0& \cdots & t^a-1 & x t^a &  \dbinom{x}{2} t^a\\
0 & \cdots & \cdots & t^a-1 & x t^a
\end{array} \right)$$
where $f$ indeed does not depend on the variable $x$. Writing $\mathbb{Q}\left[x,t\right]=\mathbb{Q}\left[x\right]\left[t\right]$ we regard $f$ and $g$ as polynomials in $t$ with coefficients in $\mathbb{Q}\left[x\right]$ and show that they do not have a common factor of positive degree. We show that $g(x,t)$ does not have a non zero complex root in $t$: let $z$ be a non zero complex number and suppose that $z^a\neq 1$; to show that $g(x,z)\in\mathbb{C}\left[x\right]$ is not the zero polynomial in $x$ we show that its term of degree one is not zero. This term is given by
\begin{eqnarray}\notag\left.\frac{\partial}{\partial x}g(x,z)\right|_{x=0} &=&\!\!
\det \left( \begin{array}{llllll}
0 & 0 &  \cdots &\cdots & 0 & \dfrac{\partial}{\partial x}\left.\left(\dbinom{x}{k}z^a\right)\right|_{x=0}\\
z^a-1 & 0 & \cdots & \cdots& 0 & \dfrac{\partial}{\partial x}\left.\left(\dbinom{x}{k-1}z^a\right)\right|_{x=0}\\
0 & z^a-1 & 0 & \cdots& 0 & \dfrac{\partial}{\partial x}\left.\left(\dbinom{x}{k-2}z^a\right)\right|_{x=0}\\
\vdots & \vdots & \ddots &\ddots & \vdots & \vdots \\
0& \cdots &  0 & z^a-1 & 0 &  \dfrac{\partial}{\partial x}\left.\left(\dbinom{x}{2}z^a\right)\right|_{x=0}\\
0 &  \cdots & \cdots & \cdots & z^a-1 & z^a
\end{array} \right)\\
\label{eq:polno0}
&= &(z^a-1)^{k-1}\dfrac{\partial}{\partial x}\left.\left(\dbinom{x}{k}z^a\right)\right|_{x=0}\neq 0
\end{eqnarray}
for every $z\in\mathbb{C}^\times$ such that $z^a\neq 1$, since $\binom{x}{k}=x(x-1)\ldots(x-k+1)k!^{-1}$ has a simple root in $x=0$. If $z^a=1$ then $g(x,z)=x^k$, which again is not the zero polynomial in $\mathbb{Q}\left[x\right]$.\\
On the other hand $t=0$ cannot be a root of $f(x,t)$, so $f$ and $g$ do not have a common root in $t$. Then their resultant $Res(f,g)$ in the variable $t$ is a non zero element $r(x)\in\mathbb{Q}\left[x\right]$ and there exist two polynomials $\phi,\psi\in\mathbb{Q}\left[x\right]\left[t\right]$ such that
$$\phi(x,t)f(x,t)+\psi(x,t)g(x,t)=r(x).$$
Therefore for every $\mu\in M_0$
\begin{eqnarray*}
&&\max\{\left|\lambda^n-1\right|_\mu , \left|\det C_{n,k}(\eta)\right|_\mu\}=
\max\{\left|f(n,\xi^n)\right|_\mu, \left|\eta^{-k}g(n,\xi^n)\right|_\mu\} \\
&\ge& \max\{\left|f(n,\xi^n)\right|_\mu, \left|g(n,\xi^n)\right|_\mu\}\\
&=& \max\{\left|f(n,\xi^n)\right|_\mu, \left|\phi(n,\xi^n)f(n,\xi^n)+\psi(n,\xi^n)g(n,\xi^n)\right|_\mu\}\\
&=& \max\{\left|f(n,\xi^n)\right|_\mu, \left|r(n)\right|_\mu\}\ge \left|r(n)\right|_\mu
\end{eqnarray*}
Then
\begin{eqnarray*}
\log\gcd(\lambda^n-1, \det C_{n,k}(\eta)) &=&
\sum_{\mu\in M_0}\log^- \max\{\left|\lambda^n-1\right|_\mu, \left|\det C_{n,k}(\eta)\right|_\mu\}\\
&\le & \sum_{\mu\in M_0}\log^-\left|r(n)\right|_\mu \le h\left(r(n)\right)=O(\log n).
\end{eqnarray*}
{\em Case 2)}
If $\eta$ is an $m$-th primitive root of unity, then
$$\det C_{n,k}(\eta)=\eta^{-k}
\det\left( \begin{array}{lllll}
n\eta^n & \dbinom{n}{2}\eta^n & \cdots & \cdots & \dbinom{n}{k}\eta^n\\
\eta^n-1 & n\eta^n & \cdots& \cdots & \dbinom{n}{k-1}\eta^n\\
\vdots & \ddots & \ddots &\vdots & \vdots \\
0& \cdots & \eta^n-1 & n\eta^n &  \dbinom{n}{2}\eta^n\\
0 & \cdots & \cdots & \eta^n-1 & n\eta^n
\end{array} \right)$$
and this is non zero for every $n\in\mathbb{N}$ sufficiently large. In fact, if $n\equiv 0\pmod{m}$, then $\det C_{n,k}(\eta)=n^k\eta^{-k}\neq 0$ for every $n\in\mathbb{N}$; otherwise, if $n\nequiv 0\pmod{m}$, we can repeat part of the above argument with minor modifications and define a polynomial $g\in\mathbb{Q}\left[x,t\right]$, with $x$ and $t$ algebraically independent over $\mathbb{Q}$, as follows:
$$g(x,t)=\det\left( \begin{array}{lllll}
x t & \dbinom{x}{2}t & \cdots & \cdots & \dbinom{x}{k}t\\
t-1 & x t & \cdots& \cdots & \dbinom{x}{k-1} t\\
\vdots & \ddots & \ddots &\vdots & \vdots \\
0& \cdots & t-1 & x t &  \dbinom{x}{2} t\\
0 & \cdots & \cdots & t-1 & x t
\end{array} \right)$$
so that
\begin{equation}
\label{eq:det}
\det C_{n,k}(\eta)=\eta^{-k} g(n,\eta^n)
\end{equation}
Let $n_0$ be an integer such that $1\le n_0\le m$, then $g(x,\eta^n)=g(x,\eta^{n_0})$ for every $n\equiv n_0\pmod{m}$. Hence, as $n$ varies, we obtain at most $m$ different polynomials $g(x,\eta),g(x,\eta^2)\ldots,g(x,\eta^m)\in\mathbb{C}\left[x\right]$ and by (\ref{eq:polno0}), if $n_0\nequiv 0\pmod{m}$, then
$$\left.\frac{\partial}{\partial x}g(x,\eta^{n_0})\right|_{x=0}=(\eta^{n_0}-1)^{k-1}\dfrac{\partial}{\partial x}\left.\left(\dbinom{x}{k}\eta^{n_0}\right)\right|_{x=0}\neq 0$$
and so $g(x,\eta^{n_0})$ is not the zero polynomial in $\mathbb{C}\left[x\right]$. Then $g(n,\eta^{n})\neq 0$ for every $n$ sufficiently large and then (\ref{eq:det}) implies $\det C_{n,k}(\eta)\neq 0$ for every $n$ sufficiently large.
Observe now that, for every $\mu\in M_0$,
\begin{equation*}
\max\{\left|\lambda^n-1\right|_\mu, \left|\det C_{n,k}(\eta)\right|_\mu\}\ge
\left|\det C_{n,k}(\eta)\right|_\mu=\left|\eta^{-k}g(n,\eta^n)\right|_\mu=
\left|g(n,\eta^n)\right|_\mu
\end{equation*}
But
\begin{equation*}
g(n,\eta^n)=\sum_{i=0}^k p_i(n)\eta^{in}
\end{equation*}
where the $p_i$ are polynomials over $\mathbb{Z}$. Then
\begin{eqnarray*}
\log\gcd(\lambda^n-1, \det C_{n,k}(\eta))&\le&
\sum_{\mu\in M_0}\log^- \left|g(n,\eta^n)\right|_\mu =
\log \prod_{\mu\in M_0}\left|g(n,\eta^n)\right|^{-1}_\mu\\
&=& \log \prod_{\mu\in M\setminus M_0}\left|g(n,\eta^n)\right|_\mu =
\log \prod_{\mu\in M\setminus M_0}\left|\sum_{i=0}^k p_i(n)\eta^{in}\right|_\mu\\
&\le & \log \prod_{\mu\in M\setminus M_0}\sum_{i=0}^k \left|p_i(n)\eta^{in}\right|_\mu = 
\log \prod_{\mu\in M\setminus M_0}\sum_{i=0}^k \left|p_i(n)\right|_\mu\\
&\le& \log \left|p(n)\right|
\end{eqnarray*}
for a suitable polynomial $p$ over $\mathbb{Z}$ and every $n$ sufficiently large. Then
$$\log\gcd(\lambda^n-1, \det C_{n,k}(\eta)) = O(\log n)$$
and this completes the proof of the lemma.
\end{proof}
We are now in a position to prove the main theorem.
\begin{proof}[Proof of Theorem \ref{mainprop}]{\em Case 1)} 
Suppose that $l-\overline{l}+h_\gamma-\overline{h}_\gamma< d-r$ for every
$\gamma\in\Gamma$. This inequality is equivalent to
$d-l-h_\gamma+\overline{l}+\overline{h}_\gamma\ge r+1$ and this in turn amounts to
say that for a chosen $\gamma\in\Gamma$, say $\gamma_1$, there exists a minor, say $y_{n,r+1,1}$, of
$B^n-I$ of order $r+1$, which is diagonal in blocks and whose blocks,
using notation of lemma \ref{lem:ckgen},
are of type $C_{n,k_i}(\lambda_i)$ where $\lambda_i\in\gamma_1\cup\Lambda^\ast$ or principal
minors of $B^n-I$ relative to eigenvalues not in $\gamma_1\cup
\Lambda^\ast$. The minor $y_{n,r+1,1}$ will thus have the form
\begin{equation}
\label{eq:y1}
y_{n,r+1,1}=\prod_{i\in\mathcal{I}}\det
C_{n,k_i}(\lambda_i)\cdot\prod_{j\in\mathcal{K}}(\eta_j^n-1),
\end{equation}
where $\mathcal{I}$ is a finite set of indexes, $\lambda_i\in\gamma_1\cup\Lambda^\ast,
\forall i\in\mathcal{I}$ and $\mathcal{K}$ is a finite set of indexes of cardinality $r+1-\sum_{i\in\mathcal{I}}k_i$ such that
$\eta_j\in\Lambda'\setminus \gamma_1$ for each $j\in\mathcal{K}$.\\
Let now $\Omega_0$ be the product of the elements of a maximal subset of cardinality at most $r+1$ of diagonal elements $\lambda^n-1$ of
$B^n-I$, where $\lambda\in\gamma_1$, i.e.
$$\Omega_0=\prod_{j\in\mathcal{L}}(\lambda_j^n-1)$$
where $\mathcal{L}$ is a finite set of indexes of cardinality at most $r+1$, $\lambda_j\in\gamma_1$ for every $j\in\mathcal{L}$, the $\lambda_j$ not necessarily distinct.
Let $\Omega_1,\ldots,\Omega_t$ be the determinants of all the minors of order $\max\{0,r+1-h_{\gamma_1}\}$
chosen from the blocks of the matrix $B^n-I$ not relative to eigenvalues in $\gamma_1$ and which do not contain elements
$\lambda^n-1$ with $\lambda\in\Lambda^\ast$; these minors exist, if $h_{\gamma_1}<r+1$, since $d-l-h_{\gamma_1}+\overline{l}\ge r+1-\overline{h}_{\gamma_1}\ge r+1-h_{\gamma_1}$. Otherwise, if $h_{\gamma_1}\ge r+1$, set $t=1$
and $\Omega_1=1$. As last, in equation (\ref{eq:y1}), set $\Omega_{t+1}\coloneqq\prod_{i\in\mathcal{I}}\det
C_{n,k_i}(\lambda_i)$ and $\Omega_{t+2}\coloneqq\prod_{j\in\mathcal{K}}(\eta_j^n-1)$.
Then
\begin{eqnarray*}
&\log\gcd_i(y_{n,r+1,i}) \le 
\log\gcd\left(\Omega_0\Omega_1,\Omega_0\Omega_2,\ldots,\Omega_0\Omega_t,
\Omega_{t+1}\Omega_{t+2}\right) \\
& \le  \log\gcd\left(\Omega_1,\Omega_2,\ldots,\Omega_t\right)+
\log\gcd\left(\Omega_0,\Omega_{t+1}\Omega_{t+2}\right)\\
& \le  \log\gcd\left(\Omega_1,\Omega_2,\ldots,\Omega_t\right)+
\log\gcd\left(\Omega_0,\Omega_{t+1}\right) +
\log\gcd\left(\Omega_0,\Omega_{t+2}\right)
\end{eqnarray*}
Observe now that
\begin{equation}
\label{eq:LI1}
\log\gcd\left(\Omega_0,\Omega_{t+1}\right) \le\sum_{j\in\mathcal{L}}
\sum_{i\in\mathcal{I}}\sum_{\mu\in M_0}\log^-
\max\{|\lambda_j^n-1|_\mu,|\det C_{n,k_i}(\lambda_i)|_\mu\}.
\end{equation}
By lemma \ref{lem:ckgen}, for every $i\in\mathcal{I}$ and $j\in\mathcal{L}$,
\begin{equation}
\label{eq:LI2}
\sum_{\mu\in M_0}\log^-
\max\{|\lambda_j^n-1|_\mu,|\det C_{n,k_i}(\lambda_i)|_\mu\}
= O(\log n),
\end{equation}
for every $n\in\mathbb{N}$ sufficiently large. Putting together equations (\ref{eq:LI1}) and (\ref{eq:LI2}) we have, for every $\epsilon>0$,
\begin{equation*}
\log\gcd\left(\Omega_0,\Omega_{t+1}\right) \le \epsilon n,
\end{equation*}
for $n$ sufficiently large.\\
Observe now that
\begin{equation*}
\log\gcd\left(\Omega_0,\Omega_{t+2}\right) \le\sum_{i\in\mathcal{L}}
\sum_{j\in\mathcal{K}}\sum_{\mu\in M_0}\log^-
\max\{|\lambda_i^n-1|_\mu,|\eta_j^n-1|_\mu\}.
\end{equation*}
Following \cite{cz:matperiod}, we can now apply the following fact, stated as Proposition 2 in \cite{cz:lowheight} (beware that our definition of $\log^-$ differes from that of \cite{cz:lowheight}, where $\log^- x=\min\{0,\log x\}$):
\begin{proposition}[Proposition 2 of \cite{cz:lowheight}]
\label{czprop}
Let $\delta>0$. All but finitely many solutions $(u,v)\in
(\mathcal{O}_{\mathbb{K},S}^*)^2$ to the inequality
$$\sum_{\mu\in M_0}\log^- \max\{|u-1|_\mu,|v-1|_\mu\}>\delta\max\{h(u),h(v)\}$$
satisfy one of finitely many relations $u^a v^b=1$, where $a,b\in\mathbb{Z}$
are not both zero.
\end{proposition}
We apply this fact with $u=\lambda_i^n$ and $v=\eta_j^n$. Since
$\lambda_i\nsim\eta_j$ for each $i\in\mathcal{L}$ and
$j\in\mathcal{K}$, then for each $\tilde{\epsilon}>0$
$$\sum_{\mu\in M_0}\log^- \max\{|\lambda_i^n-1|_\mu,|\eta_j^n-1|_\mu\}\le
\tilde{\epsilon} \max\{h(\lambda_i^n),h(\eta_j^n)\}=\tilde{\epsilon} n
\max\{h(\lambda_i),h(\eta_j)\},$$
for $n$ sufficiently large, for every $i\in\mathcal{L}$ and $j\in\mathcal{K}$. Therefore we get
\begin{equation*}
\log\gcd\left(\Omega_0,\Omega_{t+2}\right)\le \tilde{\epsilon} n \sum_{i\in\mathcal{L}}
\sum_{j\in\mathcal{K}} \max\{h(\lambda_i),h(\eta_j)\}
\end{equation*}
and taking $\tilde{\epsilon}=\epsilon\left(
\sum_{i\in\mathcal{L}} \sum_{j\in\mathcal{K}}
\max\{h(\lambda_i),h(\eta_j)\}\right)^{-1}$, we obtain
\begin{equation*}
\log\gcd\left(\Omega_0,\Omega_{t+2}\right)\le\epsilon n,
\end{equation*}
for $n$ sufficiently large.\\
If $h_{\gamma_1}\ge r+1$ the proof of case 1 can be concluded since
\begin{equation*}
\log\gcd_i(y_{n,r+1,i}) \le 
\log\gcd\left(\Omega_0,\Omega_{t+1}\right) +
\log\gcd\left(\Omega_0,\Omega_{t+2}\right) \le 
\epsilon n + \epsilon n
\end{equation*}
for $n$ sufficiently large.\\
Otherwise, if $h_{\gamma_1}< r+1$ we are left with giving a suitable upper bound  for $\log\gcd\left(\Omega_1,\Omega_2,\ldots,\Omega_t\right)$; observe that $\Omega_1,\Omega_2,\ldots,\Omega_t$ are all the minors of order $r+1-h_{\gamma_1}$ of the matrix $B^n-I$, deprived of its blocks relative to eigenvalues in $\gamma_1$, that do not contain elements
$\lambda^n-1$ with $\lambda\in\Lambda^\ast$. For every $\gamma\neq\gamma_1$ we have $d-l-h_{\gamma_1}-h_{\gamma}+\overline{l}+\overline{h}_\gamma\ge r+1-h_{\gamma_1}$, then we can repeat the procedure up to here developed, by replacing $d$ with $d-h_{\gamma_1}$, $\Gamma$ with $\Gamma\setminus \{\gamma_1\}$, $r$ with $r-h_{\gamma_1}$ and considering only the minors $\Omega_1,\Omega_2,\ldots,\Omega_t$ instead of all the minors of $B^n-I$. We come up with a new set $\{\Omega_0^1,\Omega_1^1,\ldots,\Omega_{t_1+2}^1\}$ and by possibly iterating this procedure, we come up after a finite number, say $s$, of steps with the case where $t_s=1$ and we can conclude that for every $\epsilon>0$
\begin{equation*}
\log\gcd_i(y_{n,r+1,i}) \le \epsilon n
\end{equation*}
for $n$ sufficiently large. Thus $A$ is $r$-regular.\\
{\em Case 2)} Suppose now that there exists a $\gamma\in\Gamma$, say $\gamma_1$ such that $l-\overline{l}+h_{\gamma_1}-\overline{h}_{\gamma_1} \ge d-r$. This inequality is equivalent to
$d-l-h_{\gamma_1}+\overline{l}+\overline{h}_{\gamma_1} < r+1$ and this in turn amounts to say that in the determinant of each minor of order $r+1$ of the matrix $B^n-I$ there is a factor $\lambda^n-1$ with $\lambda\in\gamma_1\cup\Lambda^\ast$.\\
Let now $T$ be the order of torsion in the subgroup of $\mathbb{K}^\ast$ generated by the eigenvalues of $B$ and observe that for each $\lambda_i\in\gamma_1$ there exist two integer $a_i,b_i$ such that $\lambda_1^{a_i}=\lambda_i^{b_i}$. Let $m$ be the least common multiple of $T$ and the $b_i$'s and consider the subset $\mathcal{N}$ of the natural numbers defined by
$$\mathcal{N}=\{n\in\mathbb{N} \text{ such that } n\equiv 0 \pmod{m}\}$$
For each $n\in\mathcal{N}$, say $n=jm$ with $j\in\mathbb{N}$ and for each $\mu\in M_0$ we have
$$\max_i\{\left|y_{n,r+1,i}\right|_\mu\}\le \left|\lambda_1^j-1\right|_\mu$$
and hence
$$\log\gcd_i(y_{n,r+1,i})\ge\sum_{\mu\in M_0}\log^-\left|\lambda_1^j-1\right|_\mu$$
for every $j\in\mathbb{N}$.
We will now prove that there exits $\rho>0$ such that
\begin{equation}
\label{eq:gtrhon}
\sum_{\mu\in M_0}\log^-\left|\lambda_1^j-1\right|_\mu>\rho j
\end{equation}
for every $j$ in an infinite subset of $\mathcal{N}$.
Observe now that
\begin{eqnarray*}
\sum_{\mu\in M_0}\log^-\left|\lambda_1^j-1\right|_\mu &=&
\sum_{\mu\in M}\log^-\left|\lambda_1^j-1\right|_\mu-
\sum_{\mu\in M\setminus M_0}\log^-\left|\lambda_1^j-1\right|_\mu \\
&=& h\left(\lambda_1^j-1\right)-\sum_{\mu\in M\setminus M_0}\log^-\left|\lambda_1^j-1\right|_\mu \\
&=& j h\left(\lambda_1\right) +O(1) -\sum_{\mu\in M\setminus M_0}\log^-\left|\lambda_1^j-1\right|_\mu
\end{eqnarray*}
Hence proving (\ref{eq:gtrhon}) amounts to prove that there exists $\rho>0$ such that
\begin{equation}
\label{eq:logm}
\sum_{\mu\in M\setminus M_0}\log^-\left|\lambda_1^j-1\right|_\mu < j\left(h\left(\lambda_1\right)-\rho \right)
\end{equation}
for every $j$ in an infinite subset of $\mathcal{N}$. The last inequality is true since we will now prove that $\forall \epsilon >0$
\begin{equation}
\label{eq:A}
\sum_{\mu\in M\setminus M_0}\log^-\left|\lambda_1^j-1\right|_\mu < \epsilon j + O(1)
\end{equation}
for every $j$ in an infinite subset of $\mathcal{N}$, by applying the (generalized) Roth's theorem \cite[chapter 6]{b:heights} in the following form.
\begin{theorem}[{\bf Roth}]
\label{thm:ss}
Let $\mathbb{K}$ be a number field and $S$ a finite set of places. For each $\mu\in S$ let $\alpha_\mu$ be $\mathbb{K}-$algebraic. Then for each $\epsilon >0$, there exist only finitely many $\beta\in\mathbb{K}$ such that
\begin{equation*}
\prod_{\mu\in S}\min\left(1,|\beta-\alpha_\mu|_\mu\right)\le H(\beta)^{-2-\epsilon}
\end{equation*}
\end{theorem}
To prove (\ref{eq:A}), let us define
\begin{eqnarray*}
D(j)&\coloneqq& \prod_{\substack{\mu\in M\setminus M_0 \\ |\lambda_1|_\mu<1}}\min\left\{1,\left|\lambda_1^j-1\right|_\mu \right\} \\
E(j)&\coloneqq&\prod_{\substack{\mu\in M\setminus M_0 \\ |\lambda_1|_\mu>1}}\min\left\{1,\left|\lambda_1^j-\infty \right|_\mu \right\} \\
F(j)&\coloneqq&\prod_{\substack{\mu\in M_0 \\ |\lambda_1|_\mu<1}}\min\left\{1,\left|\lambda_1^j-0\right|_\mu \right\},
\end{eqnarray*}
where $\left|\lambda_1^j-\infty \right|_\mu \coloneqq \left|\lambda_1^j \right|_\mu^{-1}$.
Then, for every $\epsilon_1>0$, Roth's theorem implies that
\begin{equation*}
D(j) E(j) F(j)>H\left(\lambda_1^j\right)^{-2-\epsilon_1}
\end{equation*}
for every $j$ sufficiently large. Observe now that
\begin{equation*}
E(j)=\prod_{\substack{\mu\in M\setminus M_0 \\ |\lambda_1|_\mu>1}}\frac{1}{\max\left\{1,\left|\lambda_1^j\right|_\mu \right\}}=H(\lambda_1^j)^{-1}
\end{equation*}
since $\lambda_1$ is an algebraic integer. Moreover
\begin{equation*}
F(j)=\prod_{\substack{\mu\in M_0 \\ |\lambda_1|_\mu<1}}\left|\lambda_1^j\right|_\mu=
\prod_{\mu\in M_0}\left|\lambda_1^j\right|_\mu=
\prod_{\mu\in M\setminus M_0}\left|\lambda_1^j\right|_\mu^{-1}=c^j H(\lambda_1^j)^{-1}
\end{equation*}
where
\begin{equation*}
c\coloneqq \prod_{\substack{\mu\in M\setminus M_0 \\ |\lambda_1|_\mu<1}} \left|\lambda_1\right|^{-1}
\end{equation*}
is a constant, depending on $\lambda_1$, with $c>1$. Hence, putting everything together,
\begin{equation*}
D(j)> c^{-j} H(\lambda_1)^{-j\epsilon_1}
\end{equation*}
for every $j$ sufficiently large. Let us now define $b\coloneqq c^{1/\epsilon_1}$ and observe that for every $\delta>1$ and for every $\epsilon>0$,
\begin{equation*}
c^{-j} H(\lambda_1)^{-j\epsilon_1} > \delta \exp(-j \epsilon)
\end{equation*}
for every $j$ sufficiently large, when
\begin{equation*}
\epsilon_1<\frac{\epsilon}{\log\left(b H(\lambda_1)\right)}
\end{equation*}
Hence, taking into account that
\begin{equation*}
\sum_{\mu\in M\setminus M_0}\log^-\left|\lambda_1^j-1\right|_\mu =
\log \prod_{\mu\in M \setminus M_0}\min\left\{1,\left|\lambda_1^j-1\right|_\mu \right\}^{-1}=\log D(j)^{-1}
\end{equation*}
for every $j$ sufficiently large, we conclude that
\begin{equation*}
\sum_{\mu\in M\setminus M_0}\log^-\left|\lambda_1^j-1\right|_\mu < \log\left(\delta^{-1}\exp(j\epsilon)\right)=j\epsilon + O(1)
\end{equation*}
for every $j$ sufficiently large and this proves (\ref{eq:A}).\\
We can therefore conclude that
\begin{eqnarray*}
\log\gcd_i(y_{n,r+1,i}) &\ge& \sum_{\mu\in M_0}\log^-\left|\lambda_1^j-1\right|_\mu \\
&=& j h\left(\lambda_1\right) +O(1) -\sum_{\mu\in M\setminus M_0}\log^-\left|\lambda_1^j-1\right|_\mu \\
&>& j h\left(\lambda_1\right)-j\epsilon +O(1)\\
&=&j(h(\lambda_1)-\epsilon) + O(1)> \rho j
\end{eqnarray*}
for $j$ sufficiently large, where $\rho$ is for instance $(h(\lambda_1)-\epsilon)/2$. This proves (\ref{eq:gtrhon}) and then $A$ is $r$-exceptional.
\end{proof}
\begin{proof}[Proof of Corollary \ref{cor:alggrp}]
If $A$ is $r$-exceptional, then by Theorem \ref{mainprop}, there exists $\gamma\in \Gamma$ such that
$l+h_\gamma\ge d-r+\overline{l}+\overline{h}_\gamma$.\\
If $A$ is diagonalizable, i.e. if $f=0$, then $\overline{l}+\overline{h}_\gamma=0$ and hence $l+h_\gamma\ge d-r$. Thus $e\le d-(l+h_\gamma)+1\le d-(d-r)+1=r+1$.\\
If $A$ is diagonalizable, i.e. if $f=1$, then $\overline{l}+\overline{h}_\gamma\ge 1$ and hence $l+h_\gamma\ge d-r+1$. Thus $e\le d-(l+h_\gamma)+1\le d-(d-r+1)+1=r$.\\
In both cases $A$ $r$-exceptional implies $e+f\le r+1$.
\end{proof}
Let us now come to the proof of Theorem \ref{invarprop}.
\begin{proof}[Proof of Theorem \ref{invarprop}]
Let $\phi$ be an endomorphism of a free module over a finitely generated ring $R$ of characteristic zero and let $d$ be the dimension of the module. Let $\lambda_1,\lambda_2,\ldots,\lambda_d$ be the eigenvalues of $\phi$ each repeated with its algebraic multiplicity. Finally let $\alpha_{n,1},\alpha_{n,2},\ldots,\alpha_{n,d}$ be the invariants of $\phi^n-I$, that are rational integers by hypothesis. Recalling that
\begin{eqnarray*}
\alpha_{n,1}&=&s_1(\lambda_1^n-1,\ldots,\lambda_d^n-1)\\
\alpha_{n,2}&=&s_2(\lambda_1^n-1,\ldots,\lambda_d^n-1)\\
\vdots&&\vdots\\
\alpha_{n,d}&=&s_d(\lambda_1^n-1,\ldots,\lambda_d^n-1)
\end{eqnarray*}
where $s_k$ is the $k$-th elementary symmetric polynomial, we have
\begin{equation}
\label{eq:combinv}
(\lambda_i^n-1)^d=\sum_{k=1}^d (-1)^{k+1}(\lambda_i^n-1)^{d-k}\alpha_{n,k}
\end{equation}
for every $i=1,2,\ldots,d$. Fix now a positive integer $N$ and suppose that $N|\alpha_{n,k}^{d!k^{-1}}$ for every $k=1,\ldots,d$. Then $N|\alpha_{n,k}^{d!}$ for every $k=1,\ldots,d$, and using (\ref{eq:combinv}), we have
\begin{eqnarray*}
\log\gcd_i(\lambda_i^n-1)^d & = & \sum_{\mu\in M_0}\log^-\max_i \left|(\lambda_i^n-1)^d\right|_\mu \\
& = & \sum_{\mu\in M_0}\log^-\max_i\left| \sum_{k=1}^d (-1)^{k+1}(\lambda_i^n-1)^{d-k}\alpha_{n,k} \right|_\mu \\
& \ge & \sum_{\mu\in M_0}\log^-\max_i \max_k \left| (\lambda_i^n-1)^{d-k}\alpha_{n,k} \right|_\mu \\
& \ge & \sum_{\mu\in M_0}\log^- \max_k \left|\alpha_{n,k} \right|_\mu = \log \gcd_k (\alpha_{n,k})\ge d!^{-1}\log N
\end{eqnarray*}
where $M_0$ is the set of non archimedean valuations of the field of fractions of the ring $R$.
Suppose now that $\lambda_1$ and $\lambda_2$ are two multiplicatively independent eigenvalues of $\phi$ and apply Proposition \ref{czprop}, as we did in proving Theorem \ref{mainprop}. We obtain, for every $\epsilon>0$,
\begin{equation}
\label{eq:magg2}
\log \gcd_i (\lambda_i^n-1) \le \log \gcd \left\{ \lambda_1^n-1,\lambda_2^n-1 \right\}\le \epsilon n
\end{equation}
for every $n$ sufficiently large. Therefore
\begin{equation*}
\log N \le d!d\epsilon n
\end{equation*}
for every $n$ sufficiently large, i.e.
\begin{equation*}
\lim_{N\rightarrow\infty}\frac{k(\phi,N)}{\log N}=+\infty
\end{equation*}
On the other hand if all the eigenvalues of $\phi$ are pairwise multiplicatively dependent, we can proceed as in the proof of Theorem \ref{mainprop}. For every $i=1,\ldots ,d$ there exist two integers $a_i,b_i$, not both zero, such that $\lambda_1^{a_i}=\lambda_i^{b_i}$. Let $m$ be the least common multiple of the $b_i$'s and consider the subset $\mathcal{N}$ of the natural numbers defined by
$$\mathcal{N}=\{n\in\mathbb{N} \text{ such that } n\equiv 0 \pmod{m}\}$$
For each $n\in\mathcal{N}$, say $n=jm$ with $j\in\mathbb{N}$ and for each $\mu\in M_0$ we have
$$\max_k\{\left|\alpha_{n,k}\right|_\mu\}\le \left|\lambda_1^j-1\right|_\mu$$
and applying the Roth's theorem as in the proof of Theorem \ref{mainprop} we get, for every $\epsilon>0$,
\begin{eqnarray*}
\log\gcd_k(\alpha_{n,k}^{d!k^{-1}}) &\ge&
\log\gcd_k(\alpha_{n,k}) \\
&\ge& \sum_{\mu\in M_0}\log^-\left|\lambda_1^j-1\right|_\mu 
> j(1-\epsilon) h(\lambda_1) + O(1)
\end{eqnarray*}
for $j$ sufficiently large. Hence there exists a positive constant $\rho$ such that
\begin{equation*}
\log\gcd_k(\alpha_{n,k}^{d!k^{-1}}) > \rho n
\end{equation*}
for every sufficiently large $n\in\mathcal{N}$. Then, taking $N_n\coloneqq \gcd_k(\alpha_{n,k}^{d!k^{-1}})$, we get
$$k(\phi,N_n)\le n \le \frac{1}{\rho}\log
\gcd_k(\alpha_{n,k}^{d!k^{-1}})=\frac{1}{\rho}\log N_n$$
and so
$$\frac{k(\phi,N_n)}{\log N_n}\le \frac{1}{\rho},\quad\forall n\in\mathcal{N}.$$
\end{proof}
To conclude we prove Theorems \ref{ellprop} and \ref{2ellprop}.
\begin{proof}[Proof of Theorem \ref{ellprop}]
Let $E$ be an elliptic curve over a finite field $\mathbb{F}_q$ and let $\phi:E\rightarrow E$ be the Frobenius endomorphism. Let $m(q^n), l(q^n)$ be the integers that determine the structure of the group of $\mathbb{F}_{q^n}-$rational points, as in (\ref{eq:grpstruct}), with $m(q^n)|l(q^n)$. Recall now that we may define the determinant and the trace of an endomorphism $\phi$ of an elliptic curve $E$  by choosing a prime $l$ different from the characteristic of $\mathbb{F}_q$ and considering the representation
\begin{eqnarray*}
\text{End}(E)&\rightarrow& \text{End}(T_l(E))\\
\phi &\rightarrow& \phi_l
\end{eqnarray*}
of the ring $\text{End}(E)$ of endomorphisms of $E$ into the ring of endomorphisms of the $l$-adic Tate module of $E$. Since $l$ is coprime with $q$, $T_l(E)$ is isomorphic to $\mathbb{Z}_l\times \mathbb{Z}_l$ and if we choose a basis for this $\mathbb{Z}_l$-module, we can write $\phi_l$ as a $2\times 2$ matrix whose entries belong to $\mathbb{Z}_l$. It is then possible to compute $\det(\phi_l)$ and $\text{Tr}(\phi_l)$, and it turns out that these quantities are rational integers independent from the chosen prime $l$ \cite[chapter 5]{silver:arit_ell}. We can then define
\begin{eqnarray*}
\det(\phi)&\coloneqq& \det(\phi_l) \\
\text{Tr}(\phi)&\coloneqq& \text{Tr}(\phi_l)
\end{eqnarray*}

Recall now that the Weil pairing
\begin{equation*}
e_{m(q^n)}: E[m(q^n)] \times E[m(q^n)] \rightarrow \mbox{\boldmath$\mu$}_{m(q^n)}=m(q^n)^{\text{th}}\, \text{roots of unity}
\end{equation*}
is surjective and Galois invariant (see \cite[chapter III]{silver:arit_ell}). There exists then an $m(q^n)^{\text{th}}$ primitive root of unity which belongs to the image of $e_{m(q^n)}$ and is $\mathbb{F}_{q^n}$-rational. Since the order of $\mathbb{F}_{q^n}^\ast$ is $q^n-1$, then $m(q^n)|q^n-1$ and then $m(q^n)$ is coprime with the characteristic of $\mathbb{F}_q$. In this case the subgroup $E[m(q^n)]$ of $m(q^n)-$torsion points of $E$ is isomorphic to $\mathbb{Z}/m(q^n)\mathbb{Z} \times \mathbb{Z}/m(q^n)\mathbb{Z}$ and this is a subgroup of $E(\mathbb{F}_{q^n})$, since $m(q^n)|l(q^n)$; in other words $m(q^n)$-torsion points are $\mathbb{F}_{q^n}$-rational. 

Consider now the multiplication map
\begin{eqnarray*}
[m(q^n)]&:&E\rightarrow E \\
&&P\mapsto [m(q^n)]P
\end{eqnarray*}
and the following isogeny
\begin{eqnarray*}
\phi^n-I&:&E\rightarrow E \\
&&P\mapsto \phi^n(P)-P
\end{eqnarray*}
Since $m(q^n)$ is coprime with $q$, the multiplication map $[m(q^n)]$ is separable (see \cite[chapter III]{silver:arit_ell}), and satisfies
\begin{equation*}
 \ker \left([m(q^n)]\right) \subset \ker \left(\phi^n-I\right)
\end{equation*}
There exists then a unique isogeny
\begin{equation*}
\psi_n:E\rightarrow E
\end{equation*}
such that $\phi^n-I=\psi_n\circ [m(q^n)]$ (see \cite[chapter III, Corollary 4.11]{silver:arit_ell}).
This implies that
\begin{equation}
\label{eq:mdet}
m(q^n)|\det(\phi^n-I)
\end{equation}
and
\begin{equation*}
m(q^n)|\text{Tr}\left(\phi^n-I\right)
\end{equation*}
and then
\begin{equation}
\label{eq:mtrace}
m(q^n)|\left(\text{Tr}(\phi^n-I)\right)^2
\end{equation}
If $E$ is an ordinary elliptic curve, the map $\phi$ possesses two multiplicatively independent eigenvalues $\alpha,\beta$, for otherwise $\alpha^a=\beta^b$, for suitable $(a,b)\in\mathbb{Z}^2\setminus\{(0,0)\}$, would imply $a=b$, since $|\alpha|=|\beta|=\sqrt{q}$. Then $\alpha^{2a}+\beta^{2a}\equiv 0 \pmod{q^a}$ and $E$ would be supersingular. We can then apply Theorem \ref{invarprop}: since (\ref{eq:mdet}) and (\ref{eq:mtrace}) hold, then
\begin{equation*}
\frac{n}{\log m(q^n)}\ge \frac{k(\phi,m(q^n))}{\log m(q^n)}
\end{equation*}
and therefore
\begin{equation}
\label{eq:mlim}
\lim_{n\rightarrow\infty}\frac{n}{\log m(q^n)}=+\infty
\end{equation}
But recalling the Hasse-Weil relation (\ref{eq:hasseweil})
\begin{equation*}
\sharp E(\mathbb{F}_{q^n})=l(q^n)m(q^n)=q^n+1-\text{Tr}(\phi^n)
\end{equation*}
and the fact that $\text{Tr}(\phi^n)=O(q^{n/2})$, we get
\begin{equation*}
\lim_{n\rightarrow\infty}\frac{\log l(q^n)+\log m(q^n)}{n\log q}=1
\end{equation*}
Hence (\ref{eq:mlim}) implies
\begin{equation*}
\lim_{n\rightarrow\infty}\frac{\log l(q^n)}{n\log q}=1
\end{equation*}
and this in turn implies that for every $\epsilon >0$,
\begin{equation*}
l(q^n)>q^{n(1-\epsilon)}
\end{equation*}
for every $n$ sufficiently large.\\
On the other hand if $E$ is supersingular, there exist \cite[chapter 13]{huse:ell} two strictly positive integers $a, b$ such that $\phi^a=[p^b]$, where $p=\text{char}(\mathbb{F}_q)$. Let $\mathcal{N}\coloneqq \{n\in\mathbb{N} | n\equiv 0 \pmod{a}\}$ and observe that if $n\in\mathcal{N}$, say $n=ja$, with $j\in\mathbb{N}$, then $P\in E$ is 
$\mathbb{F}_{q^n}-$rational if and only if
\begin{equation*}
0=(\phi^n-I)(P)=[p^{bj}-1](P)
\end{equation*}
i.e. if and only if $P$ is a $(p^{bj}-1)-$torsion point. But 
\begin{equation*}
E[p^{bj}-1]\cong \mathbb{Z}/(p^{bj}-1)\mathbb{Z} \times \mathbb{Z}/(p^{bj}-1)\mathbb{Z}
\end{equation*}
because $p^{bj}-1$ is coprime with $p$. Hence
\begin{equation*}
m(q^n)=p^{b\frac{n}{a}}-1
\end{equation*}
for every $n\in\mathcal{N}$. If $\forall \epsilon >0$, $l(q^n)>q^{n(1-\epsilon)}$ for every $n$ sufficiently large, then
\begin{equation}
\label{eq:ep}
\sharp E(\mathbb{F}_{q^n})=l(q^n)m(q^n)> q^{n(1-\epsilon)} \left(p^{b\frac{n}{a}}-1\right)
\end{equation}
for every $n\in\mathcal{N}$ sufficiently large and this would contradict the Hasse-Weil relation since
\begin{equation*}
\sharp E(\mathbb{F}_{q^n})<q^n+1+2q^{n/2}<2q^n
\end{equation*}
for every $n\in\mathbb{N}$ and this, together with (\ref{eq:ep}), would imply
\begin{equation*}
q^{n(1-\epsilon)} \left(p^{b\frac{n}{a}}-1\right)<2q^n
\end{equation*}
for every $n\in\mathcal{N}$ sufficiently large, leading to a contradiction when $\epsilon$ is sufficiently small.
\end{proof}
To prove Theorem \ref{2ellprop}, recall that two elliptic curves $E_1$ and $E_2$ are isogenous over $\mathbb{F}_{q}$ if and only if they have the same number of $\mathbb{F}_{q}-$rational points \cite[chapter 5]{silver:arit_ell}. So if $E_1$ and $E_2$ are two $\mathbb{F}_{q}-$isogenous elliptic curves, then the Frobenius endomorphisms $\phi_1$ and $\phi_2$ have the same characteristic polynomial and hence $\phi_1$ and $\phi_2$ have the same eigenvalues.\\
Viceversa let $\alpha_i, \overline{\alpha}_i$ be the complex conjugate eigenvalues of the Frobenius endomorphism of $E_i$, for $i=1,2$. If $\phi_1$ and $\phi_2$ have multiplicatively dependent eigenvalues, then there exists a positive integer $a$ such that $\alpha_1^a=\alpha_2^a$ and then automatically $\overline{\alpha}_1^a=\overline{\alpha}_2^a$. Hence $\phi_1^a$ and $\phi_2^a$ have the same eigenvalues and therefore the same characteristic polynomial. Then $\sharp E_1(\mathbb{F}_{q^a})=\sharp E_2(\mathbb{F}_{q^a})$ and so $E_1$ and $E_2$ are isogenous over $\mathbb{F}_{q^a}$. Observe moreover that $\alpha_1^a=\alpha_2^a$ implies $\alpha_1=\zeta \alpha_2$ for a certain $a-$th root of unity $\zeta$ which belongs to $\mathbb{Q}(\alpha_1,\alpha_2)$. Hence $\left[\mathbb{Q}(\zeta):\mathbb{Q}\right]=1,2$ or $4$. If $\left[\mathbb{Q}(\zeta):\mathbb{Q}\right]=1$ then $\zeta=\pm 1$ and $a=1$ or $2$; if $\left[\mathbb{Q}(\zeta):\mathbb{Q}\right]=2$ then $\zeta=\pm \imath, \pm \rho, \pm \rho^2$, where $\rho=\exp(2\pi \imath/3 )$ and $a=3,4$ or $6$; if $\left[\mathbb{Q}(\zeta):\mathbb{Q}\right]=4$ then $\zeta$ is a primitive root of unity of order $5,8,10$ or $12$ and consequently $a=5,8,10$ or $12$. To summarize, if $\phi_1$ and $\phi_2$ have multiplicatively dependent eigenvalues, then $E_1$ and $E_2$ are isogenous over $\mathbb{F}_{q^a}$, where $a\in\{1,2,3,4,5,6,8,10,12\}$.
\begin{proof}[Proof of Theorem \ref{2ellprop}]
 Since $E_1$ and $E_2$ are ordinary, then by Theorem \ref{ellprop} we have, $\forall \epsilon >0$,
\begin{equation}
\label{eq:ll}
l_1(q^n)l_2(q^n) \ge q^{2n(1-\epsilon)}
\end{equation}
for every $n$ sufficiently large. Let $\phi_1$ and $\phi_2$ be the Frobenius isogenies of $E_1$ and $E_2$ and let $\phi$ be the Frobenius isogeny of their product $\mathcal{A}$. We can choose a basis in $T_l(\mathcal{A})$ such that the matrix representing $\phi$ is diagonal of the form
\begin{equation}
\label{eq:fro4}
\phi_l\coloneqq
\left( \begin{array}{cccc}
\alpha_1 & 0 & 0 & 0\\
0 & \overline{\alpha}_1 & 0&  0\\
0 & 0 & \alpha_2 & 0 \\
0 & 0 & 0 & \overline{\alpha}_2
\end{array} \right)
\end{equation}
where $\alpha_i,\overline{\alpha}_i$ are the complex conjugate eigenvalues of $\phi_i$, $i=1,2$.
If $E_1$ and $E_2$ are not $\mathbb{F}_{q}-$isogenous, then by remark preceeding this proof and the fact that $E_1$ and $E_2$ are ordinary, $\alpha_1,\overline{\alpha}_1,\alpha_2,\overline{\alpha}_2$ are pairwise multiplicative independent. Hence, by Theorem \ref{mainprop} the matrix (\ref{eq:fro4}) representing $\phi_l$ is $2-$regular, in the sense that (\ref{eq:lgcd}) holds with $r=2$.
If we define $\Delta(q^n)\coloneqq\gcd(l_1(q^n),l_2(q^n))$, then $\Delta(q^n)$ divides all the determinants of the minors of order $3$ of $\phi_l^n-I$ and then $\forall \epsilon>0$
\begin{equation}
\label{eq:Delta}
\Delta(q^n)<\exp(\epsilon n)\;\;
\text{for }n\text{ sufficiently large.}
\end{equation}
We can then conlcude by (\ref{eq:ll}) and (\ref{eq:Delta}) that $\forall \epsilon>0$
\begin{equation*}
l(q^n)=\frac{l_1(q^n)l_2(q^n)}{\Delta(q^n)}>
q^{2n(1-\epsilon)}\exp(-\epsilon n)\;\;
\text{for }n\text{ sufficiently large.}
\end{equation*}
Viceversa if $E_1$ and $E_2$ are $\mathbb{F}_{q}-$isogenous, then $\alpha_1$ and $\alpha_2$ are multiplicatively dependent (possibly exchanging $\alpha_2$ with $\overline{\alpha}_2$). Then by Theorem \ref{mainprop} the matrix (\ref{eq:fro4}) is $2-$exceptional, i.e. $\exists\rho>0$ and an infinite subset $\mathcal{N}\subset \mathbb{N}$ such that, if we let $l_i'(q^n)\coloneqq l_i(q^n)/\Delta(q^n)$ for $i=1,2$, then $\forall\epsilon>0$
\begin{eqnarray*}
\rho n&<&\log\gcd\left((\alpha_1^n-1)(\overline{\alpha}_1^n-1)(\alpha_2^n-1),(\alpha_1^n-1)(\overline{\alpha}_1^n-1)(\overline{\alpha}_2^n-1),\right.\\
&&\left.(\alpha_1^n-1)(\alpha_2^n-1)(\overline{\alpha}_2^n-1),(\overline{\alpha}_1^n-1)(\alpha_2^n-1)(\overline{\alpha}_2^n-1)\right)\\
&=&\log\gcd\left(l_1(q^n)m_1(q^n)(\alpha_2^n-1),l_1(q^n)m_1(q^n)(\overline{\alpha}_2^n-1),\right.\\
&&\left.(\alpha_1^n-1)l_2(q^n)m_2(q^n),(\overline{\alpha}_1^n-1)l_2(q^n)m_2(q^n)\right)\\
&=&\log\Delta(q^n)+\log\gcd\left(l_1'(q^n)m_1(q^n)(\alpha_2^n-1),l_1'(q^n)m_1(q^n)(\overline{\alpha}_2^n-1),\right.\\
&&\left.(\alpha_1^n-1)l_2'(q^n)m_2(q^n),(\overline{\alpha}_1^n-1)l_2'(q^n)m_2(q^n)\right)\\
&\le&\log\Delta(q^n)+\log\gcd\left(\alpha_2^n-1,\overline{\alpha}_2^n-1\right)\\
&&+\log\gcd\left(l_1'(q^n)m_1(q^n),(\alpha_1^n-1)l_2'(q^n)m_2(q^n),(\overline{\alpha}_1^n-1)l_2'(q^n)m_2(q^n)\right)\\
&\le&\log\Delta(q^n)+\log\gcd\left(\alpha_2^n-1,\overline{\alpha}_2^n-1\right)\\
&&+\log\gcd\left(l_1'(q^n)m_1(q^n),l_2'(q^n)m_2(q^n)\right)+\log\gcd\left(\alpha_1^n-1,\overline{\alpha}_1^n-1\right)\\
&=&\log\Delta(q^n)+\log\gcd\left(\alpha_2^n-1,\overline{\alpha}_2^n-1\right)\\
&&+\log\gcd\left(m_1(q^n),m_2(q^n)\right)+\log\gcd\left(\alpha_1^n-1,\overline{\alpha}_1^n-1\right)\\
&\le&\log\Delta(q^n)+\epsilon n+\log m_1(q^n) + \epsilon n
\end{eqnarray*}
for every $n\in\mathcal{N}$ sufficiently large, where the last inequality follows since $\alpha_i$ and $\overline{\alpha}_i$ are multiplicatively independent, for $i=1$ and $2$. Remember now that
\begin{equation*}
m_1(q^n)<\exp(\epsilon n)
\end{equation*}
for $n$ sufficiently large, since $E_1$ is ordinary. This proves that
\begin{equation*}
\log \Delta(q^n)>\rho n - 3\epsilon n
\end{equation*}
for every $n\in\mathcal{N}$ sufficiently large. If $\rho'>0$ is a real constant, $\rho'<\rho$, then
\begin{equation*}
\Delta(q^n)>\exp(\rho' n)
\end{equation*}
for every $n\in\mathcal{N}$ sufficiently large. Hence
\begin{equation*}
l(q^n)=\frac{l_1(q^n)l_2(q^n)}{\Delta(q^n)}<
l_1(q^n)l_2(q^n)\exp(-\rho' n)
\end{equation*}
for every $n\in\mathcal{N}$ sufficiently large.\\
Moreover by the Hasse-Weil relation
\begin{equation*}
l_1(q^n)l_2(q^n)\le \sharp E_1(\mathbb{F}_{q^n})\sharp E_2(\mathbb{F}_{q^n})<(q^n+1+2q^{n/2})^2<4q^{2n}
\end{equation*}
for $n$ sufficiently large and so
\begin{equation*}
l(q^n)=<
4q^{2n}\exp(-\rho' n)
\end{equation*}
for every $n\in\mathcal{N}$ sufficiently large and this contradicts (\ref{eq:2cycl}) if
\begin{equation*}
\epsilon<\frac{1}{2}\frac{\rho'}{1+2\log q}
\end{equation*}
It is now straightforward to prove (\ref{eq:gcdnoiso}). In fact, if $E_1$ and $E_2$ are ordinary and not isogenous, then (\ref{eq:mlim}) and (\ref{eq:Delta}) imply that for every $\epsilon>0$
\begin{eqnarray*}
\gcd \left(\sharp E_1(\mathbb{F}_{q^n}),\sharp E_2(\mathbb{F}_{q^n})\right)&=&
\gcd \left(m_1(q^n)l_1(q^n),m_2(q^n)l_2(q^n)\right)\\
&\le& m_1(q^n)m_2(q^n)\Delta(q^n)\le \exp\left(\frac{\epsilon}{3}n\right)^3=\exp(\epsilon n)
\end{eqnarray*}
for every $n$ sufficiently large.
\end{proof}